\input amstex
\documentstyle{amsppt}
\def\o{operator}
\define\p{pencil}
\def\op{operator pencil}
\define\dif{differential\ }
\define\eq{equation}
\define\resp{respectively}
\def\se{\sigma_{\text{ess}}}
\def\si{\sigma}
\def\aa{\alpha}
\def\la{\lambda}
\def\bb{\beta(x)}
\def\th{\theta}
\def\lp{L(\lambda)}
\def\bn{\binom}
\redefine\Re{\operatorname{Re}}
\define\Ker{\operatorname{Ker}}
\def\sil{\sigma(L)}
\def\ga{\gamma}

\let\kap=\varkappa
\define\pd#1#2#3{\frac{\partial^{#1} #2}{\partial#3^{#1}}}
\define\od#1#2{\frac{d^{#1}}{d#2^{#1}}}
\def\H{\Cal H}
\def\D{\frak D}

\topmatter
\author Rostyslav O.~Hryniv\\
Moscow State University, Russia \endauthor
\title On the problem of semiinfinite beam oscillation with internal
damping\footnotemark"*"
    \endtitle
\footnotetext"*"{ This work was partially supported by Russian Fund of Basic
Research,
    grant No. 96-01-01292.}
\abstract
We study the Cauchy problem for the equation of the form
$$
  \ddot{u}(t) + (\aa A + B)\dot{u}(t) + (A+G)u(t) = 0
\tag*
$$
where $A$, $B$, and $G$ are \o s in a Hilbert space $\H$ with $A$
selfadjoint,
$\si(A)=[0,\infty)$, $B\geqslant0$ bounded, and $G$ symmetric and $A$-
subordinate in
a certain sense. Spectral properties of the corresponding \op\
 $\lp := \la^2I + \la (\aa A + B) + A + G$
are studied, and existence and uniqueness of generalized and classical
solutions of the Cauchy problem are proved. Equations of the type (*)
include, e.g.,
an abstract model for the problem of semiinfinite beam oscillations
with
internal damping.

This article is based on the lecture delivered at VII Crimean Autumn
Mathematical School-Sym\-po\-si\-um on Spectral and Evolutionary
Problems,
Crimea, Ukraine, September 18--29, 1996.
\endabstract
\endtopmatter

\document
\head{Introduction}\endhead
The aim of the present article is to study some class of \dif \eq s
and corresponding \op s in a Hilbert space, which provide abstract
models for
many problems in elasticity theory, hydrodynamics, control theory etc.

Consider, for example, a visco-elastic semiinfinite beam placed in
viscous
external medium. Its small transverse oscillations are described in
dimensionless coordinates by the equation (cf.~\cite{P1})
$$
  \aa\frac{\partial^5u}{\partial t \partial x^4}
    + \pd4{u}{x}
        + \pd{}{}{x}\left(g(x) \pd{}{u}{x}\right)
    + \bb\pd{}{u}{t}
    + \pd{2}{u}{t} = 0, 
    \qquad x\geqslant0,\quad t\geqslant0.
\tag1
$$
Here $u(x,t)$ is the transverse displacement of the beam at point $x$
and
time $t$;
$\aa > 0$ is a small parameter specifying internal damping,
$\bb$ determines external damping,
and $g(x)$ describes tension force distribution.

Suppose for simplicity that the left beam end is clamped, i.e. that
$u(x,t)$
satisfies the boundary conditions
$$
 u(0,t) = \pd{}{u(x,t)}{x}\mid_{x=0} =0,
\tag2
$$
and let at the moment $t=0$ the profile and velocity of the beam are
$$
  u(x,0) = \psi_0(x), \qquad \pd{}{u(x,t)}{t}\mid_{t=0} = \psi_1(x).
\tag3
$$

We will represent the initial boundary-value problem (1)--(3) as a
Cauchy
problem for an abstract equation in the Hilbert space
$L_2(0,\infty)$, and
then will study the latter by means of unbounded \o\ theory and \op\
theory.
(See also \cite{P1} and \cite{H1} for some related results.)

\head{1. Abstract \dif\ \eq\ and \op}\endhead
Problem (1)--(3) can be written in the form of
$$
\align
   &\ddot{u}(t) + (\aa A + B)\dot{u}(t) + (A+G)u(t) = 0,
\tag4\\
   &u(0) = \psi_0, \qquad \dot{u}(0) = \psi_1,
\tag5
\endalign
$$
where $u(t)$ is a function taking its values into the Hilbert space
$\H:=L_2(0,\infty)$,\, $\psi_0,\,\break \psi_1 \in \H$,
and $A,\,B$, and $G$ are linear \o s in $\H$ defined by
the equalities\footnotemark%
\footnotetext{$\D(T)$ denotes the domain of the operator $T$.}
$$
\aligned
  (Ay)(x) &= y^{\text{iv}}(x),\\(By)(x) &= \bb y(x),\\(Gy)(x) &=
            \left(g(x)y'(x)\right)',
\endaligned\qquad
\aligned
  \D(A) &= \{y \in W_2^4(0,\infty)\,|\,y(0) = y'(0) = 0\},\\
  \D(B) &= H,\\ \D(G) &= \D(A).
\endaligned
\tag6
$$

Suppose that the functions $g(x),\, g'(x)$, and $\bb$ are real,
measurable,
and essentially bounded; moreover, $g(x)$ and $\bb$ belong to ``the
class
$\Cal K$'' (\cite{B}), i.~e\. there exists a number $a > 0$ such that
$$
 \lim_{x \to \infty}\int\nolimits_{x-a}^{x+a}{|g(s)|\,ds} = 0, \qquad
 \lim_{x \to  \infty}\int\nolimits_{x-a}^{x+a}{|\beta(s)|\,ds} = 0.
$$
Then the \o s $A$, $B$ and $G$ possess the following properties:
\roster
\item"(H1)" $A=A^*>0$, the essential spectrum of the operator $A$
 coincides with the semiaxis $[0,\infty)$;
\item"(H2)" $B=B^*\geqslant0$  is bounded and
 $(A+I)$-compact in the sense of quadratic forms
 (cf. \cite{B},\cite{RS});
\item"(H3)" $G$  is symmetric, $A^{1/2}$-subordinated,
 i.e. $|G| \leqslant g_0 A^{1/2} + g_1 I$ for some $g_0,\, g_1>0$, and
 $(A+I)$-compact in the sense of quadratic forms.
\endroster

In the sequel we will study abstract Cauchy problem (4)--(5) in the
Hilbert
space $\H$ under hypotheses (H1)--(H3) on the \o s $A$, $B$ and $G$
only, and
will not exploit their concrete form (6) untill section 4.

Let there exists a solution to \eq\ (4) of the form $u(t) = e^{\la
t}y$ with
$\la \in \Bbb C$ and $y \in \H$. In order to find all such $\la$ and
$y$ we
get the spectral  problem
$$
  [\la^2I + \la (\aa A + B) + A + G] y = 0
$$
for the quadratic \op\
$$
  \lp := \la^2I + \la (\aa A + B) + A + G
$$
in the Hilbert space $\H$.

\head{2. Spectral properties of \op\ $\lp$}\endhead
Behavior of solutions to \eq\ (4) depends heavily on the structure
and
localization of \op\ $\lp$ spectrum, and so in this section
we will briefly discuss  some spectral properties of $\lp$.

First recall that the {\it spectrum} $\sil$ of the \p\ $\lp$ is the
complement
in the complex plane $\Bbb C$ to the set $\rho(L)$ of all {\it
regular} points;
here $\la_0 \in \rho(L)$ iff the \o\ $L(\la_0)$ is boundedly
invertible and
the  inverse \o\ $L^{-1}(\la_0)$ is defined on the whole space $\H$.
We
distinguish in $\sil$ the {\it point spectrum}
$$
  \si_p(L) := \left\{\la_0 \in \sil \mid\, \Ker
L(\la_0)\ne\{0\}\right\};
$$
and the {\it essential spectrum}
$$
  \se(L) := \left\{\la_0 \in \sil \mid\, \text{the \o\ } L(\la_0)
    \text{ is not a Fredholm one}\right\}.
$$
Any number $\la_0 \in \si_p(L)$ is called an {\it eigenvalue} (EV),
and
any nonzero vector $y_0 \in \Ker L(\la_0)$ is called a corresponding
{\it eigenvector} of the \p\ $\lp$.

\subhead{2.1. Essential spectrum}\endsubhead
 Let $O$ denote the circle of radius $1/\aa$
with centrum at the point $-1/\aa$ and $J$ denote the interval
$(-\infty,-1/\aa]$.

\proclaim{Theorem 1 (\cite{H1})} The essential spectrum $\se(L)$
of the pencil $\lp$ coincides with the set $O\cup{J}$.
\endproclaim

\subhead{2.2. Nonreal eigenvalues}\endsubhead
Let $\Pi_- := \{z\in \Bbb C\, |\,\Re\,z< 0\}$ be the left half-plane
and numbers\footnotemark\
\footnotetext{Some of the numbers $b_\pm$ and $g_\pm$ may equal
$\pm\infty$.}%
$b_+$ and $b_-$ ($g_+$ and $g_-$) denote the upper and lower bounds
of the
\o\ $B$ (of the \o\ $G$, \resp). Note that due to hypotheses (H1)--
(H3) we have
the inequalities $\pm b_\pm \geqslant 0$ and $\pm g_\pm \geqslant 0$; moreover,
$b_-=0$.

\proclaim{Lemma 2}
All the nonreal EV's of the pencil $\lp$ belong to the set
$\Pi_- \cap M \cap R$, where $M$ is the set\footnotemark
$$
  M:= \{\la \in \Bbb C\,\mid\,|\la + 1/\aa|^2
    \leqslant 1/\aa^2 + g_1 + \sqrt{-2\,g_0\Re\la/\aa}\}
$$
\footnotetext{The numbers $g_0$ and $g_1$ were introduced in
hypothesis (H3).}
and $R$ is the ring\footnotemark
$$
  R:= \{\la \in \Bbb C\,|\, r_- \leqslant |\la - 1/\aa|^2 \leqslant r_+\}
$$
with the numbers $r_\pm$ determined via the \o s $A$, $B$, and $G$.
In particular, if the \o\ $G$ is bounded above (below), then we can
put
$r_+ := \frac1{\aa}\sqrt{1+ \aa^2g_+}$ (\resp,
$r_-:=\frac1{\aa}\sqrt{1 - \aa b_+ + \aa^2g_-}$).
\endproclaim
\footnotetext{Which can degenerate into a disc or a point.}

\demo{Proof}
The assertion about $\Pi_-$ and $R$ (as well as the choice of $r_\pm$)
was proved in \cite{H2}.
Next, it was shown in \cite{H1} (cf. also \cite{P1}) that for any
$\ga > 1/\aa$ all the nonreal EV's are contained in the set
$M_\ga := \{\la \in \Bbb C\,\mid\,
        |\la + \ga|^2 \leqslant \ga^2 - g_1 + g_0^2/4(\ga\aa - 1)\}$,
and the intersection $\bigcap\limits_{\ga > 1/\aa} M_\ga$ is easily
seen to
coincide with the set $M$.
\enddemo
\pagebreak

\subhead{2.3. The spectrum in the right half-plane}\endsubhead
\proclaim{Lemma 3}
The nonzero spectrum of the pencil $\lp$ in the closed right half-
plane
consists of the real isolated EV's; their number $\kap_1(L)$ counted
according
to multiplicities equals\footnotemark\ $\nu(A + G)$, the total
multiplicity of
negative spectrum of the \o\ $A + G$.
\endproclaim
\footnotetext{Therefore, the {\it instability index} $\kap(L)$ of the
\p\
$\lp$, i.~e\. the number of linearly independent increasing solutions
to \eq\
(4), is not less than $\kap_1(L)$. If $\la = 0$ is a singular
critical point
(\cite{La}), then $\kap(L) > \kap_1(L)$.}

\demo{Proof}
It is a corollary of proposition 6 in \cite{LSY}. Similar result was
also proved
in \cite{LS} and \cite{P2}.
\enddemo

\proclaim{Corollary 4}
The point $\la = 0$ is an accumulation point of real EV's of the \p\
$\lp$
from the right iff $\nu(A + G) = \infty$.
\endproclaim

Note that for concrete \dif\ \o s (6) the quantity $\nu(A + G)$ can
be easily
estimated from above, see section~4.

\subhead{2.4. Accumulation of real EV's at the points $-1/\aa$ and
$0$}\endsubhead
Let for $k>-1/\aa$ a number $\nu(k)$ denote the total multiplicity of
negative spectrum of the \o\ $L(k)$.

\proclaim{Lemma 5 (\cite{Hr1})} Suppose that $\nu(-1/\aa) = \infty$
\rom($\nu(0) = \infty$\rom)\rom;
then EV's from the interval $(-1/\aa,0)$ accumulate at
the point $-1/\aa$ (at the point $0$, \resp).
\endproclaim

\head{3. The Cauchy problem}\endhead
\subhead{3.1. Classical and generalized solutions}\endsubhead
We start with the following definition.

\definition{Definition}
Let $S$ and $T$ be closed \o s in $\H$. A function
 $u(t) \in C^2(\Bbb R_+,\H)$
is said to be a {\it classical solution} to the \eq
$$
  \ddot{u}(t) + S \dot{u}(t) + Tu(t) = 0
\tag7
$$
if for any $t>0$ we have $u(t) \in \D(T)$, $\dot{u}(t) \in \D(S)$, and
equality $(7)$ holds.
\enddefinition

\define\wt{\widetilde}
\define\wL{\wt L}%
\define\wB{\wt B}
\define\wC{\wt C}
\define\wT{\wt{\Bbb T}}
\define\U{\Bbb U}
Fix a number $k_0 > \sup\limits_{\la\in \sil}\Re\la$ and consider the
\p
$$
  \wL(\xi) := L(\xi + k_0) = \xi^2 I + \xi \wB + \wC,
$$
where $\wB:= 2 k_0 I + \aa A + B \gg 0$ and $\wC:= L(k_0) \gg 0$. It
is
easily seen that a function $u(t)$ is a classical solution to \eq\
(4) iff
the function $v(t) := e^{-k_0 t} u(t)$ is a classical solution to the
\eq
$$
  \wL(\tfrac{d}{dt})v(t) := \ddot{v}(t) + \wB\dot{v}(t) + \wC v(t) =
0.
\tag8
$$
If in addition equalities (5) hold, then $v(t)$ satisfies the initial
conditions
$$
  v(0) = \psi_0, \qquad \dot{v}(0) = -k_0 \psi_0 + \psi_1.
\tag9
$$
Problem (8)--(9) now can be reduced to the first order system
$$
\align
  \dot{\pmb V}(t) &= \wt{\Bbb T}\, \pmb V(t),
\tag10\\
  \pmb V(0) &= \pmb \psi := \bn{\psi_0}{-k_0\psi_0+ \psi_1},
\tag11
\endalign
$$
in the space $\H \times \H$, where
$$
  \pmb V(t) = \bn{v_1(t)}{v_2(t)} \text{ \ and \ }
  \wt{\Bbb T} = \pmatrix  0 & I \\ - \wt C & -\wt B \endpmatrix.
$$
Actually it is more natural to consider system (10)--(11) not in the
space
$\H \times \H$, but in the so-called ``energy'' space
$\Bbb H = \H_{1/2} \times \H$, where the Hilbert space scale $\H_\th$
is
generated by the \o\ $\wC$ (namely, $\H_\th$  coincides with
$\D(\wC^{\th})$
and is equipped with the norm $\|\phi\|_\th := \|\wC^{\th}\phi\|$,
see for
details~\cite{LM}). Then the \o\ $\wT$ is closed and densely defined
on
the domain
$$
    \D(\wt{\Bbb T}) = \left\{\bn{x_1}{x_2} \in \H_{1/2} \times
\H_{1/2}
    \mid \, \wt C x_1 + \wt B x_2 \in \H_0 \right\}.
$$
Now we define a solution to \eq\ (10) to be any function
$\pmb V(t) \in C^1(\Bbb R_+, \Bbb H)$
such that $\pmb V(t) \in \D(\wT)$ for all $t>0$
and equality (10) is fulfilled.

It is easily seen that any classical solution $v(t)$ to \eq\ (8)
generates
the  solution $\pmb V(t) := \bigl(v(t), \dot{v}(t)\bigr)$ to \eq\
(10). On the
contrary, if $\pmb V(t) = \bigl( v_1(t), v_2(t)\bigr)$ is a solution
to
\eq\ (10), then the function $v_1(t)$, which formally satisfies (8),
may not be a classical solution to (8). Therefore it is natural to
call
the function $v_1(t)$ a {\it generalized solution} to \eq\ (8).

\subhead{3.2. Analyticity of the semigroup generated by the \o\
$\wT$}\endsubhead
First we will deal with generalized solution to problem (8)--(9). The
solvability of corresponding system (10)--(11) depend essentially on
the
properties of the \o\ $\wT$.

\proclaim{Theorem 6}
  The \o\ $\wT$ generates an analytic $C_0$-semigroup of contractions
$\U_t$
in the space $\Bbb H$.
\endproclaim

\demo{Proof}
 According to \cite{K}, it suffices to prove that for some constant
$C>0$ and
all $\xi \in \Bbb C$ with $\Re\xi>0$ the inequality
$$
  \|\left(\wT - \xi \Bbb I\right)^{-1}\|_{\frak B(\Bbb H)} \leqslant C/|\xi|
\tag12
$$
holds. The straightforward calculations show that the relation
$$
  \left(\wT - \xi \Bbb I\right)\bn{f_1}{f_2} = \bn{g_1}{g_2}
$$
implies%
\def\Lm{\wt L^{-1}(\xi)}%
$$
  \bn{f_1}{f_2} =  \left(\wT - \xi \Bbb I\right)^{-1} \bn{g_1}{g_2} =
  \pmatrix
  - \Lm (\wt B + \xi I) &  - \Lm \\ \Lm\wt C & -\xi \Lm
  \endpmatrix \bn{g_1}{g_2},
\tag13
$$
and henceforth (12) follows from inequalities (a)--(d) in Lemma 7
below.
The theorem is proved.
\enddemo

\def\B{\frak B}
\proclaim{Lemma 7 (\cite{H2})}
There exist positive constants $c_j,\, j=\overline{1,4}$ such that
for all
$\xi \in \Bbb C$ with $\Re\xi > 0$ the following inequalities are
satisfied:
\roster
\item "(a)" $\| \xi \Lm\|_{\B(\H_0,\H_0)} \leqslant c_1/|\xi|$;
\item "(b)" $\|\Lm \|_{\B(\H_0,\H_{1/2})} \leqslant c_2/|\xi|$;
\item "(c)" $\|\Lm \wt C\|_{\B(\H_{1/2},\H_0)} \leqslant  c_3/|\xi|$;
\item "(d)" $\|\Lm (\wt B + \xi I)\|_{\B(\H_{1/2},\H_{1/2})} \leqslant
c_4/|\xi|$.
\endroster
\endproclaim

\subhead{3.3. Solvability of the Cauchy problem}\endsubhead
Due to theorem 6 we can easily study the generalized
solutions to Cauchy problem (4)--(5).
Let $\Bbb P$ denote the orthoprojector in $\Bbb H$ onto the first
coordinate,
i.e\. $\Bbb P (x_1,x_2) = x_1$.

\proclaim{Theorem 8}  For any initial data $\psi_0 \in \H_{1/2},\,
\psi_1 \in \H$
Cauchy problem (4)--(5) has a unique generalized solution $u(t)$ such
that
$$
  \bigl( u(t),\dot{u}(t)\bigr) \to \bigl(\psi_0, \psi_1\bigr)
$$
as $t\to0$ in the norm of the space $\Bbb H$. This solution equals
$$
  u(t) = e^{k_0 t}\Bbb P \U_t \pmb \psi,
$$
where $\pmb\psi := \bigl(\psi_0, -k_0 \psi_0 + \psi_1\bigr)$,
and satisfies the inequality
$$
  \|\dot{u}(t)\|^2_\H + \bigl( \wC u(t), u(t)\bigr) \leqslant
    e^{k_0 t} \bigl(\|\psi_1\|^2_\H + (\wC \psi_0,\psi_0) \bigr).
$$
\endproclaim

It is natural that in order to get a classical solution we should
choose
``smoother'' initial
data. Indeed, let $v(t)$ be a classical solution to problem (8)--(9)
and
$v(0) = \psi_0 \in \H_1$, $\dot{v}(0) = -k_0 \psi_0 + \psi_1 \in
\H_{1/2}$.
Applying to (8) the Laplace transform and integrating by parts, we get
$$
\align
  0 &= \int_0^\infty e^{-\xi t}
      \bigl(\ddot{v}(t) + \wt B \dot{v}(t) + \wt C v(t)\bigr)\,dt =\\
    &= \wt L(\xi) \int_0^\infty e^{-\xi t} v(t)\, dt -
       \bigl((\wt B + \xi I) \psi_0 -k_0\psi_0 + \psi_1 \bigr),\quad
\Re \xi >0
\endalign
$$
Applying now the inverse Laplace transform (see \cite{V}) we arrive at
the equality
$$
  v(t) = \frac1{2\pi i} \int_{\si_0 - i\infty}^{\si_0 + i\infty}
    e^{\xi t}{\wt L}^{-1}(\xi) \bigl[(\wt B + \xi I)\psi_0 -
        k_0 \psi_0 + \psi_1 \bigr]\,d\xi, \quad \si_0 > \xi_0.
\tag14
$$
Our aim is to prove that the function $v(t)$ defined by (14)
coincides with
$\Bbb P \U_t\pmb\psi$ and is a classical solution to problem (8)--(9)
(and hence the function $u(t):= e^{k_0 t}v(t)$ is a classical
solution to
Cauchy problem (4)--(5)).

\proclaim{Theorem 9} Let $\psi_0 \in \H_1$ and $\psi_1 \in \H_{1/2}$;
then
the function $u(t):= e^{k_0 t} \Bbb P \U_t\pmb\psi$ is a classical
solution
to Cauchy problem (4)--(5).
\endproclaim

\demo{Proof}
First, the inequalities from Lemma 7 justify the possibility to apply
the
inverse Laplace transform in the form (14), as well as inclusions
$v(t) \in \D(\wC) = \H_1$ and $\dot{v}(t) \in \D(\wB)$. It remains to
prove
that the functions $v(t)$ and $\Bbb P \U_t\pmb\psi$ coincide.

Notice that the holomorphic $C_0$-semigroup $\U_t$ can be constructed
via its
generator $\wT \in \Cal A(0,\varphi)$ by means of the integral (see
\cite{K})
$$
   \U_t \pmb \phi =  \frac1{2\pi i} \int_\ga e^{\xi t}
                (\wT - \xi \Bbb I)^{-1}\pmb\psi\,d\xi,
$$
where the contour $\ga$ surrounds the sector
$S(0,\varphi):=\{\xi \in \Bbb C \mid\, |\arg \xi - \pi| \leqslant \varphi\}$
(which contains the spectrum $\si(\wT)$ of the \o\ $\wT$), and
the integral converges strongly. Therefore, due to equality (13) we
have
$$
  \Bbb P \U_t \pmb\psi = \frac1{2\pi i} \int_\ga
    e^{\xi t}{\wt L}^{-1}(\xi) \bigl[(\wt B + \xi I)\psi_0 -
        k_0 \psi_0 + \psi_1 \bigr]\,d\xi,
$$
which coincides with (14) as the integrand is an analytical function
and hence
the contour $\ga$ can be transformed into the one from (14).
The theorem is proved.
\enddemo


\head{4. Application to problem (1)--(2)}\endhead
If the \o s $A$, $B$, and $G$ are originated by system (1)--(2) and
so are
defined by (6), we can essentially refine many of the above-listed
results.
Let $g_\pm(x):=\frac12\bigl(|g(x)| \pm g(x)\bigr)$.

\proclaim{Lemma 10} (cf\. Lemma 2) \rom{(a)} If $g_-(x) \equiv 0$,
then all the
nonreal \rom{EV's} of the pencil $\lp$ belong to the disc
$$
  D:= \{ \la \in \Bbb C \mid\, |\la + 1/\aa| \leqslant 1/\aa\}.
$$

\noindent \rom{(b)} If $g_+(x)\equiv 0$, then the nonreal spectrum of
the \p\
$\lp$ lies outside of the disc
$$
  D':= \{ \la \in \Bbb C \mid\, |\la + 1/\aa| < \sqrt{1-\aa
b_+}/\aa\},
$$
where $b_+ := \text{ess}\sup \bb$.
\endproclaim

According to Lemma 3, the number $\kap_1(L)$ of (real) EV's in the
right
half-plane equals $\nu(T)$, the total multiplicity of negative
spectrum of
the \o\ $T:= A + G = \od{4}{x} + \frac{d}{dx}(x)\frac{d}{dx}$ with the
domain $\D(T) = \{ y \in W^4_2 (\Bbb R_+)\mid\, y(0)=y'(0) =0\}$.
Consider
also the Schr\"odinger \o\ $S:= - \od{2}{x} -g(x)$ with the domain
$\D(S) = \{ y \in W^2_2(\Bbb R_+)\mid\, y(0) = 0\}$, and by $\nu(S)$
denote
the (possibly infinite) number of its negative EV's. The crucial role
in
estimating of $\nu(T)$  plays the following statement.

\proclaim{Lemma 11 (\cite{H1})} The following inequality holds:
$$
  \nu(T) \leqslant \nu(S) \leqslant \nu(T) + 1.
$$
\endproclaim

\proclaim{Corollary 12} \rom{(a)} Suppose that $g(x)\geqslant 0$ for all
sufficiently large $x$ and
$$
    \sup_{x\geqslant0}\,t\int_t^{\infty}g_+(s)\,ds =\infty
$$
Then $\nu(T) = \infty$, and $\la=0$ is an accumulation point of real
EV's
from both sides.

\noindent \rom{(b)} If $\max\limits_{x\geqslant
a}\,x\int_x^{\infty}g_+(s)\,ds \leqslant 1/4$
for some $a>0$, then $\, \nu(T)<\infty$, and EV's do not accumulate
at the
point $0$ from the right. In particular, we then have
$$
    \nu(T) \leqslant \int_0^{\infty}xg_+(x)\,dx.
$$
\endproclaim

\demo{Proof} Analogous statements for the Schr\"odinger \o\ $S$ are
well-known,
see, e.g., \cite{B} and \cite{RS}.
\enddemo

What concerns accumulation of the real EV's at the point $\la=-
1/\aa$, we have
the following result.

\proclaim{Lemma 13}
Suppose that $g_+(x) \not\equiv0$; then $\la =-1/\aa$ is an
accumulation point
of real EV's of the pencil $\lp$ from the right. If $b_+ < 1/\aa$,
then this
condition is a necessary one for accumulation.
\endproclaim

Finally, we can describe the properties of solution to problem (1)--
(3)
in terms of initial data.

\proclaim{Theorem 14}
Suppose that
$$
 \psi_0(x) \in W_{2,U}^2(\Bbb R_+) := \{ y(x) \in W_2^2(\Bbb R_+)
    \mid\, y(0) = y'(0) = 0\}
$$
and  $\psi_1(x) \in L_2(\Bbb R_+)$. Then problem (1)--(3) has a unique
generalized solution $u(x,t)$. If, in addition, $\psi_0(x)$
belongs to $W_2^4(\Bbb R_+)$, and  $\psi_1(x)$ belongs to
$W_{2,U}^2(\Bbb R)$,
then the solution $u(x,t)$ is a classical one.
\endproclaim

\subhead{Acknowledgements}\endsubhead
The author is deeply grateful to Prof. A.~A.~Shkalikov for
very useful and stimulating discussions.

\bigskip
{\smc Department of Mechanics and Mathematics, Moscow State
University,
    Moscow, 119899 Russia

rhryniv\@nw.math.msu.su }

\Refs

\widestnumber \key{LSY}

\ref\key B\by Birman M.\,S.
\paper On the spectrum of singular boundary problems
\jour Mat. USSR Sbornik
\yr 1961
\vol 55(97)
\issue 2
\pages 125--174
\lang Russian
\endref

\ref\key  H1 \by Hryniv R.\,O.
\paper On the spectrum of \op\ arising in the problem of semiinfinite
beam
    oscillations with internal damping
\jour Moscow Univ. Math Bulletin, Ser 1.
\yr 1996
\vol
\issue 1
\pages 19--23
\lang Russian
\endref

\ref\key  H2  \by Hryniv R.\,O.
\book Spectral analysis of \op\ arising in hydrodynamic problems
\bookinfo PhD thesis
\yr 1996
\publaddr Moscow
\publ Moscow State University
\lang Russian
\endref

\ref\key K\by Kato T.
\book Perturbation Theory for Linear Operators, \rom{2nd edition}
\publaddr New York
\publ Springer-Verlag
\yr 1976
\endref

\ref\key  LS \by Lancaster~P. and Shkalikov~A.\,A.
\paper Damped vibrations of beams and related spectral problems
\jour Canad. Appl. Math. Quart.
\yr 1994
\vol 2
\issue 1
\pages 45--90
\endref

\ref\key LSY\by Lancaster P., Shkalikov A.\,A., and Ye Q.
\paper Strongly definitizable linear pencils in Hilbert space
\jour Integr. Equat. Oper. Th.
\yr 1993
\vol 17
\pages 338--360
\endref


\ref\key  L \by Langer H.
\paper Spectral functions of definitizable \o s in Krein spaces
\inbook Lecture Notes in Mathematics
\yr 1982
\vol 948
\issue
\pages 1--46
\publaddr Berlin
\publ Springer-Verlag
\endref

\ref\key LM \by Lions J.-L. and  Magenes E.
\book Probl\`emes aux limites non homog\`enes et applications
\rom{Vol. 1}
\publ Paris
\publaddr Dunod
\yr 1968
\transl English transl.
\publ Springer-Verlag
\yr 1972
\endref


\ref\key P1\by  Pivovarchik V.\,N.
\paper On the vibration of a semiinfinite beam with internal and
external damping
\jour Prikl. Mat. Mekh.
\yr 1988
\vol 52
\issue 5
\pages 829--836
\lang Russian
\endref

\ref\key  P2 \by Pivovarchik V.\,N.
\paper On positive spectra of one class of polynomial \op s
\jour Integr. Equat. Oper. Th.
\yr 1994
\vol 19
\issue
\pages 314--326
\endref


\ref\key RS \by Reed M. and Simon B.
\book Methods of Modern Mathematical Physics, vol. IV
\publaddr London
\publ Academic Press
\yr 1978
\endref


\ref\key  V \by Vladimirov V.\,S.
\book Mathematical physics equation
\publaddr Moscow
\publ Nauka
\yr 1988
\lang Russian
\endref



\endRefs
\enddocument